\documentstyle[12pt]{article}  
\def\sq{\hbox {\rlap{$\sqcap$}$\sqcup$}}
\def\sq{\hbox {\rlap{$\sqcap$}$\sqcup$}}
\def\R{ {\rm R \kern -.31cm I \kern .15cm}}
\def\C{ {\rm C \kern -.15cm \vrule width.5pt \kern .12cm}}
\def\Z{ {\rm Z \kern -.27cm \angle \kern .02cm}}
\def\N{ {\rm N \kern -.26cm \vrule width.4pt \kern .10cm}}
\def\1{{\rm 1\mskip-4.5mu l} }
\def\lsim{\raise0.3ex\hbox{$<$\kern-0.75em\raise-1.1ex\hbox{$\sim$}}}
\def\gsim{\raise0.3ex\hbox{$>$\kern-0.75em\raise-1.1ex\hbox{$\sim$}}}
\def\noi{\noindent}

\def\beq{\begin{equation}}   \def\eeq{\end{equation}}
\def\bea{\begin{eqnarray}}  \def\eea{\end{eqnarray}}
\def\nn{\nonumber}
\def\noi{\noindent}
\def\beeq{\begin{eqnarray}} \def\eeeq{\end{eqnarray}}
\newcommand\mysection{\setcounter{equation}{0}\section}

\newcounter{hran}

\begin{document} 
\centerline{\large\bf Long range scattering for the Wave-Schr\"odinger system} 
 \vskip 3 truemm \centerline{\large\bf with large wave data and small Schr\"odinger data\footnote{Work supported in part by NATO Collaborative Linkage Grant 979341}} 

\vskip 0.5 truecm

\centerline{\bf J. Ginibre}
\centerline{Laboratoire de Physique Th\'eorique\footnote{Unit\'e Mixte de
Recherche (CNRS) UMR 8627}}  \centerline{Universit\'e de Paris XI, B\^atiment
210, F-91405 ORSAY Cedex, France}
\vskip 3 truemm
\centerline{\bf G. Velo}
\centerline{Dipartimento di Fisica, Universit\`a di Bologna}  \centerline{and INFN, Sezione di
Bologna, Italy}

\vskip 1 truecm

\begin{abstract}
We study the theory of scattering for the Wave-Schr\"odinger system
with Yukawa type coupling in space dimension 3. We prove in particular
the existence of modified wave operators for that system with no size
restriction on the wave data in the framework of a direct method which
requires smallness of the Schr\"odinger data, and we determine the
asymptotic behaviour in time of solutions in the range of the wave
operators.\end{abstract}

\vskip 3 truecm
\noi AMS Classification : Primary 35P25. Secondary 35B40, 35Q40, 81U99.  \par \vskip 2 truemm

\noi Key words : Long range scattering, modified wave operators, Wave-Schr\"odinger system.\par 
\vskip 1 truecm

\noindent LPT Orsay 04-36\par
\noindent May 2004\par

\newpage
\pagestyle{plain}
\baselineskip 18pt

\mysection{Introduction}
\hspace*{\parindent} This paper is devoted to the theory of scattering
and more precisely to the construction of modified wave operators for
the Wave-Schr\"odinger system (WS)$_3$ in space dimension 3, namely
\beq \label{1.1e}
\left \{ \begin{array}{l}  i\partial_t u = - (1/2)
\Delta u + A u \\ \\ \sq A = - |u|^2   \end{array} \right . 
 \eeq

\noi where $u$ and $A$ are respectively a complex valued and a real
valued function defined in space time ${I\hskip-1truemm R}^{3+1}$,
$\Delta$ is the Laplacian in ${I\hskip-1truemm R}^3$ and $\sq =
\partial_t^2 - \Delta$ is the d'Alembertian in ${I\hskip-1truemm
R}^{3+1}$. That system is Lagrangian and admits a number of formally
conserved quantities, among which the $L^2$ norm of $u$ and the energy 
\beq
\label{1.2e}
E(u,A) = \int dx \left \{ (1/2) \left ( |\nabla u|^2 + (\partial_t A)^2 + |\nabla A|^2\right ) + A|u|^2 \right \} \ .
\eeq

\noi The Cauchy problem for the (WS)$_3$ system is known to be globally
well posed in the energy space $X_e = H^1 \oplus \dot{H}^1 \oplus L^2$
for $(u, A, \partial_tA)$ \cite{1r}. \par

A large amount of work has been devoted to the theory of scattering for nonlinear equations
and systems centering on the Schr\"odinger equation, in particular for nonlinear Schr\"odinger
(NLS) equations, Hartree equations, Klein-Gordon-Schr\"odinger (KGS), Wave-Schr\"odinger (WS) 
and Maxwell-Schr\"odinger (MS) systems. As in the case of the linear Schr\"odinger
equation, one must distinguish the short range case from the long range case. In the former
case, ordinary wave operators are expected and in a number of cases proved to exist, describing
solutions where the Schr\"odinger function behaves asymptotically like a solution of the free
Schr\"odinger equation. In the latter case, ordinary wave operators do not exist and have to be
replaced by modified wave operators including a suitable phase in their definition. In that
respect, the (WS)$_3$ system (\ref{1.1e}) belongs to the borderline
(Coulomb) long range case, because of the $t^{-1}$ decay in
$L^{\infty}$ norm of solutions of the wave equation. Such is the case
also for the Hartree equation with $|x|^{-1}$ potential. Both are
simplified models for the more complicated Maxwell-Schr\"odinger system (MS)$_3$
in ${I\hskip-1truemm R}^{3+1}$, which belongs to the same case, as well
as the Klein-Gordon-Schr\"odinger system (KGS)$_2$ in ${I\hskip-1truemm R}^{2+1}$.
\par

The construction of modified wave operators for the previous long range
equations and systems has been tackled by two methods. The first one
was initiated in \cite{10r} on the example of the NLS equation in
${I\hskip-1truemm R}^{1+1}$ and subsequently applied to the NLS
equation in ${I\hskip-1truemm R}^{2+1}$ and ${I\hskip-1truemm R}^{3+1}$
and to the Hartree equation \cite{2r}, to the (KGS)$_2$ system
\cite{11r} \cite{12r} \cite{13r} \cite{14r}, to the (WS)$_3$ system
\cite{15r} and to the (MS)$_3$ system \cite{16r} \cite{18r}. That
method is rather direct, starting from the original equation or system.
It will be sketched below. It is restricted to the (Coulomb) limiting
long range case, and requires a smallness condition on the asymptotic
state of the Schr\"odinger function. Early applications of the method
required in addition a support condition on the 
Fourier transform of the Schr\"odinger asymptotic state and a smallness
condition of the Klein-Gordon or Maxwell field in the case of the
(KGS)$_2$ or (MS)$_3$ system respectively \cite{11r} \cite{18r}. The
support condition was subsequently removed for the (KGS)$_2$ and
(MS)$_3$ system and the method was applied to the (WS)$_3$ system
without a support condition, at the expense of adding a correction term
to the Schr\"odinger asymptotic function \cite{12r} \cite{15r}
\cite{16r}. Finally the smallness condition of the KG field was removed
for the (KGS)$_2$ system, first with and then without a support
condition \cite{13r} \cite{14r}. All the previous papers on (KGS)$_2$,
(WS)$_3$ and (MS)$_3$ use spaces of fairly regular solutions, with at
least $H^2$ regularity for the Schr\"odinger function. \par

In the present paper, we reconsider the same problem for the (WS)$_3$
system in the framework of the previous method. Our purpose is twofold.
First we show that no smallness condition is required on the wave
field. Second, we treat the problem in function spaces that are as
large as possible, namely with regularity as low as possible, and with
convergence in time as slow as possible. In particular we treat the
problem with regularity of the Schr\"odinger function at the level of
$L^2$ only. Furthermore, only a weak convergence in time of the
solutions to their asymptotic form is needed, namely $t^{-\lambda}$
with $\lambda > 3/8$. Under such a weak condition, neither a support
condition nor a correction term for the asymptotic Schr\"odinger
function is needed, as long as $\lambda \leq 1/2$. We also treat the
same problem at the level of $H^1$ and $H^2$ for the Schr\"odinger
function. This does not require any reinforcement of the smallness
condition of the Schr\"odinger asymptotic state or of the time decay.
\par

In a subsequent paper, we shall treat the same problem for the (MS)$_3$
system in the Coulomb gauge in the framework of the present method.
Again no smallness condition will be required for the Maxwell field,
and a weak time decay $t^{-\lambda}$ with $\lambda > 3/8$ will be
sufficient, so that no support condition or correction term will be
needed. (On that problem see \cite{16r} \cite{18r}).\par

For completeness and although we shall not make use of that fact in the
present paper, we mention that the same problem for the Hartree
equation and for the (WS)$_3$ and (MS)$_3$ system can also be treated
by a more complex method where one first applies a phase-amplitude
separation to the Schr\"odinger function. The main interest of that
method is to remove the smallness condition on the Schr\"odinger
function, and to go beyond the Coulomb limiting case for the Hartree
equation. That method has been applied in particular to the (WS)$_3$ system and to
the (MS)$_3$ system in a special case \cite{5r} \cite{6r} \cite{7r}.
\par

We now sketch briefly the method of construction of the modified wave
operators initiated in \cite{10r}. That construction basically consists
in solving the Cauchy problem for the system (\ref{1.1e}) with infinite
initial time, namely in constructing solutions $(u,A)$ with prescribed
asymptotic behaviour at infinity in time. We restrict our attention to
time going to $+\infty$. That asymptotic behaviour is imposed in the
form of suitable approximate solutions $(u_a,A_a)$ of the system
(\ref{1.1e}). The approximate solutions are parametrized by data $(u_+,
A_+, \dot{A}_+)$ which play the role of (actually would be in simpler
e.g. short range cases) initial data at time zero for a simpler
evolution. One then looks for exact solutions $(u,A)$ of the system
(\ref{1.1e}), the difference of which with the given asymptotic ones
tends to zero at infinity in time in a suitable sense, more precisely,
in suitable norms. The wave operator is then defined traditionally as
the map $\Omega_+ : (u_+, A_+, \dot{A}_+) \to (u,A,\partial_t A)(0)$.
However what really matters is the solution $(u, A)$ in the
neighborhood of infinity in time, namely in some interval $[T, \infty
)$, and we shall restrict our attention to the construction of such
solutions. Continuing such solutions down to $t = 0$ is a somewhat
different question, connected with the global Cauchy problem at finite
times, which we shall not touch here. That problem is well controlled
for the (WS)$_3$ system, but not for instance for the (MS)$_3$ system.
\par

The construction of solutions $(u,A)$ with prescribed asymptotic
behaviour $(u_a,A_a)$ is performed in two steps. \\

\noi \underbar{Step 1}. One looks for $(u,A)$ in the form $(u,A) = (u_a+v, A_a + B)$. The system satisfied by $(v,B)$ is
\beq \label{1.3e}
\left \{ \begin{array}{l}  i\partial_t v = - (1/2)
\Delta v + A v + Bu_a - R_1 \\ \\ \sq B = - (|v|^2 + 2 \ {\rm Re} \ \overline{u}_a v) - R_2   \end{array} \right . 
 \eeq

\noi where the remainders $R_1$, $R_2$ are defined by
\beq \label{1.4e}
\left \{ \begin{array}{l}  R_1= i\partial_t u_a + (1/2)
\Delta u_a - A_a u_a \\ \\ R_2 = \sq A_a + |u_a|^2   \ . \end{array} \right . 
 \eeq

\noi It is technically useful to consider also the partly linearized system for functions $(v',B')$
\beq \label{1.5e}
\left \{ \begin{array}{l}  i\partial_t v' = - (1/2)
\Delta v' + A v' + Bu_a - R_1 \\ \\ \sq B' = - (|v|^2 + 2 \ {\rm Re} \ \overline{u}_a v) - R_2  \ . \end{array} \right . 
 \eeq

\noi The first step of the method consists in solving the system
(\ref{1.3e}) for $(v, B)$, with $(v, B)$ tending to zero at infinity in
time in suitable norms, under assumptions on $(u_a, A_a)$ of a general
nature, the most important of which being decay assumptions on the
remainders $R_1$ and $R_2$. That can be done as follows. One first
solves the linearized system (\ref{1.5e}) for $(v',B')$ with given $(v,
B)$ and initial data $(v', B')(t_0) = 0$ for some large finite $t_0$.
One then takes the limit $t_0 \to \infty$ of that solution, thereby
obtaining a solution $(v',B')$ of (\ref{1.5e}) which tends to zero at
infinity in time. That construction defines a map $\phi : (v,B) \to
(v',B')$. One then shows by a contraction method that the map $\phi$
has a fixed point. That first step will be performed in Section 2. \\

\noi \underbar{Step 2.} The second step of the method consists in
constructing approximate asymptotic solutions $(u_a, A_a)$ satisfying
the general estimates needed to perform Step~1. With the weak time
decay allowed by our treatment of Step 1, one can take the simplest
version of the asymptotic form used in previous works \cite{5r}
\cite{6r} \cite{15r}. Thus we choose \beq
\label{1.6e}
u_a = MD \exp (- i \varphi ) w_+
\eeq 

\noi where
\beq
\label{1.7e}
M \equiv M(t) = \exp (ix^2/2t) \ ,
\eeq
\beq
\label{1.8e} 
D(t) = (it)^{-n/2} D_0(t) \quad , \quad ( D_0 (t) f) (x) = f(x/t) \ ,
\eeq

\noi $\varphi$ is a real phase to be chosen below and $w_+ = Fu_+$. We
furthermore choose $A_a$ so that $R_2 = 0$, namely $A_a = A_0 + A_1$
where $A_0$ is the solution of the free wave equation $\sq A_0 = 0$
given by 
\beq
\label{1.9e}
A_0 = \cos \omega t \ A_+ + \omega^{-1} \sin \omega t \ \dot{A}_+
\eeq

\noi where $\omega = (- \Delta)^{1/2}$, and where 
\beq
\label{1.10e}
A_1(t) = \int_t^{\infty} dt' \omega^{-1} \sin (\omega (t-t'))|u_a (t')|^2 \ .
\eeq

\noi Substituting (\ref{1.6e}) into (\ref{1.10e}) yields
\beq
\label{1.11e}
A_1(t) = t^{-1} D_0(t) \widetilde{A}_1
\eeq

\noi where
\beq
\label{1.12e}
\widetilde{A}_1 = - \int_1^{\infty} d\nu \ \nu^{-3} \omega^{-1} \sin (\omega (\nu - 1)) D_0 (\nu ) |w_+|^2 \ .
\eeq

\noi In particular $\widetilde{A}_1$ is constant in time. We finally
choose $\varphi$ by imposing $\partial_t \varphi = t^{-1}
\widetilde{A}_1$, $\varphi (1) = 0$, namely
\beq
\label{1.13e}
\varphi = (\ell n\ t) \widetilde{A}_1 \ .
\eeq

We shall show in Section 3 that the previous choice fulfills the
conditions needed for Step 1, under suitable assumptions on the
asymptotic state $(u_+, A_+, \dot{A}_+)$.\par

In order to state our results we introduce some notation. We denote by
$F$ the Fourier transform and by $\parallel \cdot \parallel_r$ the norm
in $L^r \equiv L^r ({I\hskip-1truemm R}^3)$, $1 \leq r \leq \infty$.
For any nonnegative integer $k$ and for $1 \leq r \leq \infty$, we
denote by $W_r^k$ the Sobolev spaces 
$$W_r^k = \left \{ u : \parallel u; W_r^k\parallel \ = \sum_{\alpha : 0 \leq |\alpha | \leq k} \parallel \partial_x^{\alpha} u \parallel_r \ < \infty \right \}$$

\noi where $\alpha$ is a multiindex, so that $H^k = W_2^k$. We shall
need the weighted Sobolev spaces $H^{k,s}$ defined for $k$, $s \in
{I\hskip-1truemm R}$ by 
$$H^{k,s} = \left \{ u : \parallel u; H^{k,s}\parallel \ = \ \parallel (1 + x^2)^{s/2} (1 - \Delta )^{k/2} u \parallel_2 \ < \infty \right \}$$

\noi so that $H^k = H^{k,0}$. For any interval $I$, for any Banach
space $X$ and for any $q$, $1 \leq q \leq \infty$, we denote by $L^q(I,
X)$ (resp. $L_{loc}^q(I,X)$) the space of $L^q$ integrable (resp.
locally $L^q$ integrable) functions from $I$ to $X$ if $q < \infty$ and
the space of measurable essentially bounded (resp. locally essentially
bounded) functions from $I$ to $X$ if $q = \infty$. For any $h \in
{\cal C} ([1, \infty ), {I\hskip-1truemm R}^+)$, non increasing and
tending to zero at infinity and for any interval $I \subset [1, \infty
)$, we define the spaces 
\bea
\label{1.14e}
&&X(I) = \Big \{ (v, B):v\in {\cal C}(I, L^2), \parallel (v, B);X(I)\parallel \ \equiv \ \mathrel{\mathop {\rm Sup}_{t \in I }}\ h(t)^{-1} \nn \\
&&\left ( \parallel v(t)\parallel_2 \ + \ \parallel v; L^{8/3}(J,L^4)\parallel\ + \ \parallel B; L^4(J, L^4)\parallel \right ) < \infty \Big \} \ ,
\eea
\bea
\label{1.15e}
&&X_1(I) = \Big \{ (v, B):v\in {\cal C}(I, H^1), \ \nabla B, \partial_t B \in {\cal C}(I, L^2),\nn\\
&&\parallel (v, B);X_1(I)\parallel \ \equiv \ \mathrel{\mathop {\rm Sup}_{t \in I }}\ h(t)^{-1} \Big ( \parallel v(t);H^1\parallel \ + \ \parallel v; L^{8/3}(J,W_4^1)\parallel\nn \\
&& + \ \parallel B; L^4(J, L^4)\parallel\ + \ \parallel \nabla B(t) \parallel_2 \ + \ \parallel \partial_t B(t) \parallel_2 \Big ) < \infty \Big \} \ ,
\eea

\bea
\label{1.16e}
&&X_2(I) = \Big \{ (v, B):v\in {\cal C}(I, H^2) \cap {\cal C}^1(I, L^2), \ \nabla B, \partial_t B \in {\cal C}(I, L^2),\nn\\
&&\parallel (v, B);X_2(I)\parallel\ \equiv \ \mathrel{\mathop {\rm Sup}_{t \in I }}\ h(t)^{-1} \Big ( \parallel v(t);H^2\parallel \ + \ \parallel \partial_t v(t) \parallel_2 \ + \ \parallel v; L^{8/3}(J,L^4)\parallel \nn \\
&&+ \ \parallel \partial_t v;  L^{8/3}(J,L^4)\parallel \ + \ \parallel B; L^4(J, L^4)\parallel\ + \ \parallel \nabla B(t) \parallel_2 \ + \ \parallel \partial_t B(t) \parallel_2 \Big ) < \infty \Big \} \nn \\
\eea

\noi where $J = [t, \infty ) \cap I$. \par

We can now state our results.\\

\noi{\bf Proposition 1.1.} {\it Let $h(t) = t^{-1/2}$. Let $u_a$ be
defined by (\ref{1.6e}) with $w_+ = Fu_+ \in L^4$ and $c_4 = \parallel
w_+\parallel_4$ sufficently small and with $\varphi$ defined by
(\ref{1.12e}) (\ref{1.13e}). Let $A_a$ be defined by $A_a = A_0 + A_1$
with $A_0$ and $A_1$ defined by (\ref{1.9e}) and
(\ref{1.10e})-(\ref{1.12e}).\par

(1) Let $u_+ \in H^{0,2}$, let $A_+$, $\omega^{-1} \dot{A}_+ \in L^2$
and $\nabla^2 A_+$, $\nabla \dot{A}_+ \in L^1$. Then there exists $T$,
$1 \leq T < \infty$ and there exists a unique solution $(u,A)$ of the
system (\ref{1.1e}) such that $(u - u_a, A-A_a) \in X([T, \infty
))$.\par

(2) Let $u_+ \in H^{0,3} \cap H^{1,2}$, let $A_+$, $\omega^{-1}
\dot{A}_+ \in H^1$ and $\nabla^2A_+$, $\nabla \dot{A}_+ \in W_1^1$.
Then there exists $T$, $1 \leq T < \infty$, and there exists a unique
solution $(u, A)$ of the system (\ref{1.1e}) such that $(u - u_a, A-
A_a) \in X_1 ([T, \infty ))$. Furthermore $A$ satisfies the estimate
\beq
\label{1.17e}
\parallel \nabla (A-A_a)(t)\parallel_2\ \vee\ \parallel \partial_t (A - A_a)(t) \parallel_2\ \leq\ C\ t^{-3/4}
\eeq

\noi for some constant $C$ and for all $t \geq T$.\par

(3) Let $u_+ \in H^{1,3} \cap H^{2,2}$, let $A_+$, $\omega^{-1}
\dot{A}_+ \in H^1$ and $\nabla^2A_+$, $\nabla \dot{A}_+ \in W_1^1$.
Then there exists $T$, $1 \leq T < \infty$ and there exists a unique
solution $(u, A)$ of the system (\ref{1.1e}) such that $(u - u_a, A-
A_a) \in X_2 ([T, \infty ))$. Furthermore $u - u_a \in L^{8/3}([T,
\infty ), W_4^2)$ and $(u,A)$ satisfies the estimates (\ref{1.17e}) and
\beq
\label{1.18e}
\parallel \Delta (u - u_a); L^{8/3} ([t, \infty ), L^4 ) \parallel \ \leq C\ t^{-1/2}
\eeq

\noi for some constant $C$ and for all $t \geq T$.}\\

\noi {\bf Remark 1.1} The only smallness condition bears on $c_4$ and
appears at the level of the $L^2$ theory in Part (1) of Proposition
1.1. In particular there is no smallness condition on $A_a$.
Furthermore no additional smallness condition is required for the
theories at the level of $H^1$ and $H^2$.

\mysection{The Cauchy problem at infinite initial time}
\hspace*{\parindent} In this section we perform the first step of the
construction of solutions of the system (\ref{1.1e}) as described in
the introduction, namely we construct solutions $(v, B)$ of the system
(\ref{1.3e}) defined in a neighborhood of infinity in time and tending
to zero at infinity under suitable regularity and decay assumptions on
the asymptotic functions $(u_a, A_a)$ and on the remainders $R_i$. As
mentioned in the introduction, we offer three theories with $u$ (or
$v$) at the level of regularity of $L^2$, $H^1$ and $H^2$ respectively.
As a preliminary to that study, we need to solve the Cauchy problem
with finite initial time for the linearized system (\ref{1.5e}). That
system consists of two independent equations. The second one is simply
a wave equation with an inhomogeneous term and the Cauchy problem with
finite or infinite initial time for it is readily solved under suitable
assumptions on the inhomogeneous term, which will be fulfilled in the
applications. The first one is a Schr\"odinger equation with time
dependent real potential and time dependent inhomogeneity which we
rewrite in a more concise form and with slightly different notation as
\beq \label{2.1e} i\partial_t v = - (1/2) \Delta v + Vv + f \ . 
\eeq

We first give some preliminary results on the Cauchy problem with
finite initial time for that equation at the level of regularity of
$L^2$, $H^1$ and $H^2$. Those results rely in an essential way on the
well known Strichartz inequalities for the Schr\"odinger equation
\cite{3r} \cite{9r} \cite{19r} which we recall for completeness. We
define 
\beq \label{2.2e} U(t) = \exp (i(t/2) \Delta ) \ . 
\eeq

\noi A pair of H\"older exponents $(q, r)$ will be called admissible if
$0 \leq 2/q\ = 3/2 - 3/r \leq 1$. For any $r$, $1 \leq r \leq \infty$, we
define $\overline{r}$ by $1/r + 1/\overline{r} = 1$.\\

\noi {\bf Lemma 2.1.} {\it The following inequalities hold. \par

(1) For any admissible pair $(q, r)$ and for any $u \in L^2$ 
\beq
\label{2.3e} \parallel U(t) u ; L^q ({I\hskip-1truemm R}, L^r)
\parallel \ \leq \ C \parallel u \parallel_2 \ . 
\eeq

(2) Let $I$ be an interval and let $t_0 \in I$. Then for any admissible
pairs $(q_i, r_i)$, $i = 1,2$,} 
\beq \label{2.4e} \parallel 
\int_{t_0}^t dt' \ U(\cdot - t') f(t') ; L^{q_1}(I,L^{r_1})\parallel \
\leq \ C \parallel f;L^{\overline{q}_2}(I, L^{\overline{r}_2})\parallel
\ . 
\eeq

The basic result on the Cauchy problem for (\ref{2.1e}) can be stated
as follows. \\

\noi {\bf Proposition 2.1.} {\it Let $I$ be an interval and let $t_0
\in I$.\par

(1) Let $V \in L_{loc}^1(I, L^{\infty}) + L_{loc}^p(I,L^s) +
L_{loc}^{\infty}(I,L^{3/2 + \varepsilon })$ for some $p$, $1 \leq p <
\infty$ with $2/p = 2 - 3/s$ and for some $\varepsilon > 0$, and let $f
\in L^1(I, L^2) + L^2(I, L^{6/5})$. Let $v_0 \in L^2$. Then there
exists a unique solution $v$ of (\ref{2.1e}) with $v(t_0) = v_0$ such
that $v \in L_{loc}^q(I,L^r)$ for all admissible pairs $(q, r)$.
Furthermore $v \in {\cal C}(I, L^2)$ and for all $t \in I$, $v$
satisfies the equality 
\beq \label{2.5e} \parallel v(t)\parallel_2^2\ -
\ \parallel  v_0\parallel_2^2\ = \int_{t_0}^t dt'\ 2\ {\rm Im}\
<v(t'),f(t')> \ . 
\eeq

(2) Let $V$ and $f$ satisfy the assumptions of Part (1) and in addition
$\nabla V \in L_{loc}^1(I, L^3 + L^{\infty}) + L_{loc}^4(I, L^{6/5})$
and $\nabla f \in L_{loc}^1(I, L^2) + L_{loc}^2(I, L^{6/5})$. Let $v_0
\in H^1$. Then the solution $v$ of (\ref{2.1e}) obtained in Part (1)
satisfies in addition $\nabla v \in L_{loc}^q(I, L^r)$ for all
admissible pairs $(q,r)$. Furthermore $v \in {\cal C}(I,H^1)$ and for
all $t\in I$, $v$ satisfies the equality 
\beq \label{2.6e} \parallel
\nabla v(t)\parallel_2^2\ - \ \parallel  \nabla v_0\parallel_2^2\ =
\int_{t_0}^t dt'\ 2\ {\rm Im}\ <\nabla v(t'),(\nabla V)(t') v(t') +
\nabla f(t')> \ . 
\eeq

(3) Let $V$ satisfy $V \in {\cal C}(I,L^2 + L^{\infty})$, $\partial_t V
\in L_{loc}^1(I, L^2 + L^{\infty}) + L_{loc}^2(I, L^{6/5})$, and let $f$
satisfy $f \in {\cal C} (I,L^2)$, $\partial_t f \in L_{loc}^1(I, L^2) +
L_{loc}^2(I, L^{6/5})$. Let $v_0 \in H^2$. Then the solution $v$ of
(\ref{2.1e}) obtained in Part (1) satisfies in addition $\partial_t v
\in L_{loc}^q(I, L^r)$ for all admissible pairs $(q, r)$. Furthermore
$v \in {\cal C}(I,H^2) \cap {\cal C}^1(I,L^2)$ and for all $t \in I$,
$v$ satisfies the equality 
\bea \label{2.7e} 
&&\parallel \partial_t
v(t)\parallel_2^2\ - \ \parallel  - (1/2) \Delta v_0 + V(t_0) v_0 +
f(t_0)\parallel_2^2\nn\\ &&= \int_{t_0}^t dt'\ 2 \ {\rm Im} \
<\partial_t v(t') , (\partial_t V) (t') v(t') + \partial_t f(t')> \ .
\eea

\noi If in addition $V \in L_{loc}^{q_0}(I, L^{r_0} + L^{\infty})$ and
$f \in L_{loc}^{q_0}(I, L^{r_0})$ for some admissible pair $(q_0,
r_0)$, then $\Delta v \in L_{loc}^q(I,L^r)$ for all admissible pairs
$(q,r)$ with $2 \leq r \leq r_0$.}\\

The proof is a variation of that given in \cite{8r} \cite{19r}, using
extensively Lemma 2.1. For any interval $J$, let 
\beq \label{2.8e} Z(J)
= {\cal C}(J, L^2) \cap L^2(I, L^6) \ . 
\eeq

\noi The local Cauchy problem for (\ref{2.1e}) is treated by a
contraction method applied to the integral equation associated with
(\ref{2.1e}). The relevant spaces for the contraction have $v \in Z(J)$
for Part (1), $v$, $\nabla v \in Z(J)$ for Part (2) and $v$,
$\partial_t v \in Z(J)$, $v\in {\cal C}(J,H^2)$ for Part (3), for
suitable small $J$. The extension of local solutions to global ones is
easy because the problem is linear.\par

We now begin the construction of solutions of the system (\ref{1.3e}).
For any $T$, $t_0$ with $1 \leq T < t_0 \leq \infty$, we denote by $I$
the interval $I = [T, t_0]$ and for any $t \in I$, we denote by $J$ the
interval $J = [t, t_0]$. In all this section, we denote by $h$ a
function in ${\cal C}([1, \infty ), {I\hskip-1truemm R}^+)$ such that
for some $\lambda > 0$, the function $\overline{h}(t) \equiv
t^{\lambda}h(t)$ is non increasing and tends to zero as $t \to \infty$.
\par

We shall make repeated use of the following lemma. \\

\noi {\bf Lemma 2.2.} {\it Let $1 \leq q$, $q_k \leq \infty$ ($1 \leq k \leq n$) be such that 
$$\mu \equiv 1/q - \sum_k 1/q_k \geq 0 \ .$$

\noi Let $f_k \in L^{q_k}(I)$ satisfy
\beq
\label{2.9e}
\parallel f_k;L^{q_k}(J) \parallel\ \leq N_k\ h(t)
\eeq

\noi for $1 \leq k \leq n$, for some constants $N_k$ and for all $t \in I$. \par

Let $\rho \geq 0$ be such that $n\lambda + \rho > \mu$. Then the following inequality holds for all $t\in I$}
\beq
\label{2.10e}
\parallel \left ( \prod_k f_k \right ) t^{-\rho} ; L^q(J) \parallel\ \leq C \left ( \prod_k N_k\right ) h(t)^n \ t^{\mu - \rho}
\eeq

\noi where
\beq
\label{2.11e}
C = \left ( 1 - 2^{-q(n\lambda + \rho - \mu )}\right )^{-1/q} \ .
\eeq
\vskip 5 truemm

\noi {\bf Proof.} For $t \in I$, we define $I_j = [t2^j, t2^{j+1}] \cap
I$ so that $J = \ \displaystyle{\mathrel{\mathop {\cup}_{j\geq 0}}}\
I_j$. We then rewrite $L^q(J) = \ell_j^q(L^q(I_j))$. We estimate
$$\parallel \left ( \prod_k f_k \right ) t^{-\rho} ; L^q(J) \parallel\ \leq\ \parallel \left ( \prod_k \parallel f_k ;L^{q_k}(I_j)\parallel\right ) \parallel t^{-\rho} ; L^{1/\mu}(I_j) \parallel \ ; \ell_j^q \parallel$$
$$\leq \ \left ( \prod_k N_k \right ) \parallel h(t2^j)^n (t2^j)^{-\rho + \mu}; \ell_j^q \parallel$$
$$\leq \ \left ( \prod_k N_k \right ) \overline{h}(t)^n \ t^{-n\lambda - \rho + \mu} \parallel 2^{j(-n\lambda - \rho + \mu )};\ell_j^q \parallel$$
 
\noi from which (\ref{2.10e}) follows. \par \nobreak \hfill $\sq$ \\

\noi {\bf Remark 2.1.} In some special cases, the dyadic decomposition
is not needed for the proof of Lemma 2.1. For instance if all the $q_k$
are infinite, one can estimate
\bea
\label{2.12}
&&\parallel h(t)^n \ t^{-\rho} \parallel_q \ \leq \overline{h}(t)^n \parallel t^{-\rho - n \lambda} \parallel_q\nn \\
&&\leq C\ \overline{h}(t)^n\ t^{-\rho - n \lambda + 1/q} = C\ h(t)^n\ t^{-\rho + \mu}
\eea

\noi by a direct application of H\"older's inequality in $J$. The same situation occurs if $\rho > \mu$.\\

In addition to the Strichartz inequalities for the Schr\"odinger
equation (Lemma 2.1), we shall need special cases of the Strichartz
inequalities for the wave equation \cite{4r} \cite{9r}. Let $I$ be an
interval, let $t_0 \in I$ and let $B(t_0) = \partial_t B(t_0) = 0$. Then
\beq
\label{2.13e}
\parallel B;L^4(I,L^4) \parallel \ \leq \ C \parallel  \sq B; L^{4/3} (I, L^{4/3})\parallel  \ ,
\eeq

\beq
\label{2.14e}
\mathrel{\mathop {\rm Sup}_{t\in I}}\ \left ( \parallel \nabla B(t) \parallel_2 \ \vee \ \parallel \partial_t B(t) \parallel_2 \right ) \ \leq \ \parallel \sq B; L^1(I,L^2)\parallel \ .
\eeq

We now construct solutions of the system (\ref{1.3e}) at the level of
regularity of $L^2$ for $v$. The result can be stated as follows.\\

\noi {\bf Proposition 2.2.} {\it Let $h$ be defined as above with
$\lambda = 3/8$ and let $X(\cdot )$ be defined by (\ref{1.14e}). Let
$u_a \in L^{\infty}([1, \infty ), L^4)$, $A_a \in L^{\infty}([1, \infty
), L^{\infty})$, $R_1 \in L^{\infty}([1, \infty ), L^2)$ and $R_2 \in
L^{4/3}([1, \infty ) L^{4/3})$ satisfy the estimates 
\beq
\label{2.15e}
\parallel u_a(t) \parallel_4\ \leq c_4 \ t^{-3/4}\ ,
\eeq
\beq
\label{2.16e}
\parallel A_a(t) \parallel_{\infty} \ \leq a\ t^{-1}\ ,
\eeq
\beq
\label{2.17e}
\parallel R_1;L^1([t, \infty ), L^2) \parallel\ \leq r_1\ h(t)\ ,
\eeq
\beq
\label{2.18e}
\parallel R_2;L^{4/3}([t, \infty ), L^{4/3}) \parallel\ \leq r_2\ h(t)\ ,
\eeq

\noi for some constants $c_4$, $a$, $r_1$, $r_2$ with $c_4$
sufficiently small and for all $t \geq 1$. Then there exists $T$, $1
\leq T < \infty$ and there exists a unique solution $(v, B)$ of the
system (\ref{1.3e}) in the space $X([T, \infty ))$.} \\

\noi {\bf Proof.} We follow the sketch given in the introduction. Let
$1 \leq T < \infty$ and let $(v,B)\in X([T, \infty ))$. In particular
$(v, B)$ satisfies
\beq
\label{2.19e}
\parallel v(t) \parallel_2 \ \leq N_0\ h(t)
\eeq
\beq
\label{2.20e}
\parallel v ; L^{8/3}([t, \infty), L^4) \parallel \ \leq N_1\ h(t)
\eeq
\beq
\label{2.21e}
\parallel B ; L^{4}([t, \infty), L^4) \parallel \ \leq N_2\ h(t)
\eeq

\noi for some constants $N_i$ and for all $t \geq T$. We first
construct a solution $(v',B')$ of the system (\ref{1.5e}) in $X([T,
\infty ))$. For that purpose, we take $t_0$, $T < t_0 < \infty$ and we
solve the system (\ref{1.5e}) in $X(I)$ where $I = [T, t_0]$ with
initial condition $(v', B')(t_0) = 0$. Let $(v'_{t_0}, B'_{t_0})$ be
the solution thereby obtained. The existence of $v'_{t_0}$ follows from
Proposition 2.1, part (1) with $V = A$ and $f = Bu_a - R_1$. We want to
take the limit of $(v'_{t_0}, B'_{t_0})$ as $t_0 \to \infty$ and for
that purpose we need estimates of $(v'_{t_0}, B'_{t_0})$ in $X(I)$ that
are uniform in $t_0$. Omitting the subscript $t_0$ for brevity we
define 
\beq
\label{2.22e}
N'_0 =  \ \mathrel{\mathop {\rm Sup}_{t\in I}}\ h(t)^{-1} \parallel v'(t) \parallel_2
\eeq 
\beq
\label{2.23e}
N'_1 =  \ \mathrel{\mathop {\rm Sup}_{t\in I}}\ h(t)^{-1} \parallel v';L^{8/3}(J,L^4) \parallel
\eeq 
\beq
\label{2.24e}
N'_2 =  \ \mathrel{\mathop {\rm Sup}_{t\in I}}\ h(t)^{-1} \parallel B';L^{4}(J,L^4) \parallel
\eeq 

\noi where $J = [t, \infty ) \cap I$. We first estimate $N'_0$. From (\ref{2.5e}) we obtain 
\bea
\label{2.25e}
\parallel v'(t) \parallel_2 & \leq & \parallel B u_a - R_1 ; L^1(J, L^2) \parallel \nn \\
&\leq & \parallel\ \parallel B \parallel_4 \ \parallel u_a\parallel_4 \ + \ \parallel R_1 \parallel_2\ ;L^1(J) \parallel\nn \\
&\leq & C_0 \left ( c_4\ N_2 + r_1 \right ) h(t)
\eea

\noi by Lemma 2.2, so that 
\beq
\label{2.26e}
N'_0 \leq C_0 \left ( c_4\ N_2 + r_1 \right ) \ .
\eeq

We next estimate $N'_1$. By Lemma 2.1 
\bea \label{2.27e} &&\parallel v' ; L^{8/3}(J, L^4) \parallel \ \leq \ C
\left ( \parallel A_a v';L^1(J,L^2)\parallel \ + \ \parallel
Bv';L^{8/5}(J,L^{4/3})\parallel \right . \nn \\ 
&&+ \ \parallel Bu_a - R_1; L^1(J,L^2)
\parallel  \ . \eea

\noi The last norm has already been estimated by (\ref{2.25e}) while
$$\parallel A_a v';L^1(J,L^2)\parallel\ \leq \ \parallel \ \parallel A_a \parallel_{\infty} \ \parallel v'\parallel_2 \ ; L^1(J) \parallel\ \leq C\ a\ N'_0\ h(t)\ , $$
\beq
\label{2.28e}
\parallel Bv';L^{8/5}(J,L^{4/3})\parallel\ \leq \ \parallel \ \parallel B \parallel_{4} \ \parallel v'\parallel_2 \ ; L^{8/5}(J) \parallel\ \leq C\ N_2 \ N'_0\ \overline{h}(t)\ h(t)
\eeq

\noi by Lemma 2.2. Substituting (\ref{2.28e}) into (\ref{2.27e}) and using (\ref{2.26e}), we obtain
\beq
\label{2.29e}
N'_1 \leq C_1 \left ( c_4 \ N_2 + r_1\right ) \left ( 1 + a + N_2\ \overline{h}(T)\right ) \ .
\eeq

\noi We finally estimate $N'_2$. From (\ref{2.13e}), we obtain
\bea
\label{2.30e}
&&\parallel B';L^4(J,L^4)\parallel\ \leq \ C\parallel |v|^2 + 2\ {\rm Re}\ \overline{u}_av + R_2;L^{4/3}(J,L^{4/3})\parallel\nn \\
&&\leq \ C \left ( \parallel \ \parallel v \parallel_2 \left ( \parallel v \parallel_4 \ + \ \parallel u_a \parallel_4\right );L^{4/3}(J) \parallel\ + r_2\ h(t)\right )\nn \\
&&\leq C_2 \left ( c_4\ N_0 + r_2 + N_0\ N_1 \ \overline{h}(t) \right ) h(t)
\eea

\noi by Lemma 2.2, so that 
\beq
\label{2.31e}
N'_2 \leq C_2 \left (  c_4\ N_0 + r_2 + N_0\ N_1 \ \overline{h}(T) \right )  \ .
\eeq

It follows from (\ref{2.26e}) (\ref{2.29e}) and (\ref{2.31e}) that
$(v', B')$ is bounded in $X(I)$ uniformly in $t_0$.\par

We now take the limit $t_0 \to \infty$ of $(v'_{t_0}, B'_{t_0})$,
restoring the subscript $t_0$ for that part of the argument. Let $T <
t_0 < t_1 < \infty$ and let $(v'_{t_0}, B'_{t_0})$ and $(v'_{t_1},
B'_{t_1})$ be the corresponding solutions of (\ref{1.5e}). From the
$L^2$ norm conservation of the difference $v'_{t_0} - v'_{t_1}$ and
from (\ref{2.25e}), it follows that for all $t \in [T, t_0]$
\beq
\label{2.32e}
\parallel v'_{t_0}(t) - v'_{t_1}(t) \parallel_2  \ = \ \parallel v'_{t_1} (t_0)  \parallel_2\ \leq \ C_0 \left ( c_4\ N_2 + r_1 \right ) h(t_0)
\eeq

\noi while from (\ref{2.13e}) (\ref{2.30e}) and the initial conditions
\bea
\label{2.33e}
&&\parallel B'_{t_0} - B'_{t_1};L^{4}([T, t_0],L^{4})\parallel\ \leq \ C\parallel |v|^2 + 2 \ {\rm Re} \ \overline{u}_a v + R_2 ;L^{4/3}([t_0, t_1], L^{4/3}) \parallel \nn \\
&&\leq C_2 \left ( c_4\ N_0 + r_2 + N_0\ N_1 \ \overline{h}(T)\right )  h(t_0) \ .
\eea

\noi It follows from (\ref{2.32e}) (\ref{2.33e}) that there exists
$(v',B') \in L_{loc}^{\infty}([T, \infty ) ,L^2) \oplus L_{loc}^4([T,
\infty ) ,$ $L^4)$ such that $(v'_{t_0}, B'_{t_0})$ converges to $(v',
B')$ in that space when $t_0 \to \infty$. From the uniformity in $t_0$
of the estimates (\ref{2.25e}) (\ref{2.30e}), it follows that $(v',
B')$ satisfies the same estimates in $[T, \infty )$ namely that $(v', B')$ satisfies (\ref{2.26e}) (\ref{2.31e}) with $N'_i$
defined by (\ref{2.22e})-(\ref{2.24e}) with $I = [T, \infty )$.
Furthermore it follows from (\ref{2.29e}) by a standard compactness
argument that $(v', B') \in X([T, \infty ))$ and that $v'$ also
satisfies (\ref{2.29e}). Clearly $(v', B')$ satisfies the system
(\ref{1.5e}). \par

From now on, $I$ denotes the interval $[T, \infty )$. The previous
construction defines a map $\phi : (v, B) \to (v', B')$ from $X(I)$ to
itself. The next step consists in proving that the map $\phi$ is a
contraction on a suitable closed bounded set ${\cal R}$ of $X(I)$. We
define ${\cal R}$ by the conditions (\ref{2.19e})-(\ref{2.21e}) for
some constants $N_i$ and for all $t \in I$. We first show that for a
suitable choice of $N_i$ and for sufficiently large $T$, the map $\phi$
maps ${\cal R}$ into ${\cal R}$. By (\ref{2.26e}) (\ref{2.29e})
(\ref{2.31e}) it suffices for that purpose that
\beq
\label{2.34e}
\left \{ \begin{array}{l} 
(N'_0 \leq ) \ C_0 \left ( c_4 \ N_2 + r_1 \right ) \leq N_0 \\ \\ (N'_1 \leq )   \ C_1 \left ( c_4\ N_2 + r_1 \right ) \left ( 1 + a + N_2 \ \overline{h}(T)\right ) \leq N_1 \\ \\ (N'_2 \leq ) \ C_2 \left ( c_4\ N_0 + r_2 + N_0 \ N_1\ \overline{h}(T)\right ) \leq N_2 \ .
\end{array} \right .
\eeq
\noi We fulfill those conditions by choosing the $N_i$ according to
\beq
\label{2.35e}
\left \{ \begin{array}{l} 
N_0 = C_0 \left ( c_4 \ N_2 + r_1 \right ) \\ \\ N_1 =  C_1 \left ( c_4\ N_2 + r_1 \right ) (2+a)  \\ \\ N_2 =  C_2 \left ( c_4\ N_0 + r_2 + 1 \right )  
\end{array} \right .
\eeq

\noi which is possible under the smallness condition $C_0C_2 c_4^2 < 1$, and by taking $T$ sufficiently large so that 
\beq
\label{2.36e}
N_2\ \overline{h}(T) \leq 1 \quad , \quad N_0\ N_1\ \overline{h}(T) \leq 1 \ .
\eeq

\noi We next show that the map $\phi$ is a contraction on ${\cal R}$.
let $(v_i, B_i) \in {\cal R}$, $i = 1,2$, and let $(v'_i, B'_i) = \phi
((v_i, B_i))$. For any pair of functions $(f_1, f_2)$ we define $f_{\pm}
= (1/2)(f_1 \pm f_2)$ so that $(fg)_{\pm} = f_+ g_{\pm} + f_-
g_{\mp}$. In particular $u_+ = u_a + v_+$, $u_- = v_-$, $A_+ = A_a +
B_+$ and $A_- = B_-$. Corresponding to (\ref{1.5e}), $(v'_-, B'_-)$
satisfies the system
\beq
\label{2.37e}
\left \{ \begin{array}{l} i\partial_t v'_- = - (1/2) \Delta v'_- + A_+ v'_- + B_- u_a + B_- v'_+\\ \\ \sq B'_- = - 2\ {\rm Re}\ \left ( \overline{u}_a + \overline{v}_+ \right ) v_-\ .\end{array} \right .
\eeq

\noi Since ${\cal R}$ is convex and stable under $\phi$, $(v_+, B_+)$
and $(v'_+, B'_+)$ belong to ${\cal R}$, namely satisfy
(\ref{2.19e})-(\ref{2.21e}). Let $N_{i^-}$ and $N'_{i^-}$ be the seminorms of $(v_-, B_-)$ and $(v'_-, B'_-)$ corresponding to
(\ref{2.22e})-(\ref{2.24e}), namely the constants obtained by replacing
$(v', B', N'_i)$ by $(v_-, B_-, N_{i^-})$ and $(v'_-, B'_-, N'_{i^-})$
in (\ref{2.22e})-(\ref{2.24e}). We have to estimate the $N'_{i^-}$ in
terms of the $N_{i^-}$. The estimates are essentially the same as those
of $N'_i$ in terms of $N_i$ with minor differences~: the contribution
of the remainders disappear, the linear terms are the same, and the
quadratic terms are in general obtained by polarization. The only
exceptions to that rule are the $B_- v'_+$ term in the estimate of
$N'_{0^-}$ and the $\overline{v}_+v_-$ term in the estimate of
$N'_{2^-}$. Those terms are estimated as follows
\bea
\label{2.38e}
&&\parallel B_-v'_+;L^1(J,L^2)\parallel \ \leq \ \parallel\ \parallel B_-\parallel_4 \ \parallel v'_+\parallel_4\ ;L^1(J)\parallel \nn \\
&&\leq C\ N_{2^-}\ N_1\ \overline{h}(t)\ h(t)\ ,
\eea
\bea
\label{2.39e}
&&\parallel \overline{v}_+v_-;L^{4/3}(J,L^{4/3})\parallel \ \leq \ \parallel\ \parallel v_+\parallel_4 \ \parallel v_-\parallel_2\ ;L^{4/3}(J)\parallel \nn \\
&&\leq C\ N_1\ N_{0^-}\ \overline{h}(t)\ h(t)\ .
\eea

\noi We finally obtain 
\beq 
\label{2.40e}
\left \{ \begin{array}{ll} N'_{0^-} &\leq C_0 \left ( c_4 + N_1 \
\overline{h}(T)\right ) N_{2^-}\\ \\ N'_{1^-} &\leq C \left ( \left (
c_4 + N_1 \ \overline{h}(T)\right ) N_{2^-} + \left ( a + N_2 \
\overline{h}(T)\right ) N'_{0^-} \right )\\ \\ &\leq C_1 \left ( c_4 + N_1 \
\overline{h}(T)\right ) \left ( 1 + a + N_2 \ \overline{h}(T) \right )
N_{2^-}\\ \\ N'_2 &\leq C_2 \left ( c_4 + N_1 \ \overline{h}(T)\right )
N_{0^-} \ . \end{array} \right . \eeq

It follows from (\ref{2.40e}) that the map $\phi$ is a contraction for
the pair of semi norms $N_{0^-}$, $N_{2^-}$ under the condition $C_0C_2
(c_4 + N_1 \overline{h}(T))^2 < 1$ which is the combination of a
smallness condition for $c_4$ together with a condition that $T$ be
sufficiently large. The semi norm $N_{1^-}$ does not take part in the
contraction, but is controlled separately by the previous ones. The
constants $C_0$ and $C_2$ appearing in (\ref{2.40e}) can be taken to be
the same as in (\ref{2.34e}) because they are determined by the linear
terms, which are the same in both cases. There might occur additional
constants coming from the nonlinear terms. They have been omitted. This
completes the proof of the existence part of the Proposition.
Uniqueness follows from (\ref{2.40e}) with $N'_{i^-} = N_{i^-}$. \par
\nobreak \hfill $\sq$ \par

We now turn to the construction of solutions of the system (\ref{1.3e})
at the level of regularity of $H^1$ for $v$. The result can be stated
as follows. \\

\noi {\bf Proposition 2.3.} {\it Let $h$ be defined as previously with $\lambda = 3/8$ and let $X_1(\cdot )$ be defined by (\ref{1.15e}). Let $u_a$, $A_a$, $R_1$, $R_2$ satisfy the conditions (\ref{2.15e})-(\ref{2.18e}) and in addition
\beq
\label{2.41e}
\parallel u_a(t) \parallel_{\infty} \ \leq c\ t^{-3/2} \quad , \quad \parallel \nabla u_a (t) \parallel_4 \ \leq c\ t^{-3/4} \ ,
\eeq
\beq
\label{2.42e}
\parallel \nabla A_a \parallel_{\infty} \ \leq a\  t^{-1}\ ,
\eeq
\beq
\label{2.43e}
\parallel \nabla R_1 ; L^1([t, \infty ), L^2)\parallel \ \leq r_1\ h(t) \ ,
\eeq
\beq
\label{2.44e}
\parallel R_2;L^1([t, \infty ), L^2 ) \parallel \ \leq r_2\ t^{-1/2}\ h(t)
\eeq

\noi for some constants $c_4$, $c$, $a$, $r_1$, $r_2$ with $c_4$
sufficiently small and for all $t \geq 1$. Then there exists $T$, $1
\leq T < \infty$ and there exists a unique solution $(v, B)$ of the
system (\ref{1.3e}) in the space $X_1([T, \infty ))$. Furthermore $B$ satisfies the estimate
$$\parallel \nabla B(t) \parallel_2\ \vee\ \parallel \partial_t B(t)\parallel_2\ \leq C \left ( t^{-1/2} + t^{1/4} h(t)\right ) h(t)\eqno(2.44e)$$

\noi for some constant $C$ and for all $t \geq T$.}\\

\noi {\bf Proof.} The proof follows closely that of Proposition 2.2. Let
$1 \leq T < \infty$ and let $(v, B) \in X_1([T, \infty ))$. In
particular $(v, B)$ satisfies (\ref{2.19e})-(\ref{2.21e}) and in
addition
\beq
\label{2.45e}
\parallel \nabla v(t) \parallel_{2} \ \leq \  N_3\ h(t) 
\eeq
\beq
\label{2.46e}
\parallel \nabla v ; L^{8/3}([t, \infty ), L^4)\parallel \ \leq \  N_4\ h(t) 
\eeq
\beq
\label{2.47e}
\parallel \nabla B(t)\parallel_2 \ \vee \ \parallel \partial_t B(t)\parallel_2\ \leq \  N_5\ h(t) 
\eeq

\noi for some constants $N_i$ and for all $t \geq T$. We first
construct a solution $(v',B')$ of the system (\ref{1.5e}) in $X_1([T,
\infty ))$. For that purpose, we take $t_0$, $T < t_0 < \infty$ and we
solve the system (\ref{1.5e}) in $X_1(I)$ where $I = [T, t_0]$ with
initial condition $(v', B')(t_0) = 0$. Let $(v'_{t_0}, B'_{t_0})$ be
the solution thereby obtained. The existence of $v'_{t_0}$ follows from
Proposition 2.1, part (2) with $V = A$, and $f = B u_a - R_1$. We want
to take the limit of $(v'_{t_0}, B'_{t_0})$ as $t_0 \to \infty$ and for
that purpose we need estimates of $(v'_{t_0}, B'_{t_0})$ in $X_1(I)$
that are uniform in $t_0$. Omitting the subscript $t_0$ for brevity we
define $N'_i$, $0 \leq i \leq 5$, by (\ref{2.22e})-(\ref{2.24e}) and by 
\beq
\label{2.48e}
N'_3 =  \ \mathrel{\mathop {\rm Sup}_{t\in I}}\ h(t)^{-1} \parallel \nabla v'(t) \parallel_2
\eeq 
\beq
\label{2.49e}
N'_4 =  \ \mathrel{\mathop {\rm Sup}_{t\in I}}\ h(t)^{-1} \parallel \nabla v';L^{8/3}(J,L^4) \parallel
\eeq 
\beq
\label{2.50e}
N'_5 =  \ \mathrel{\mathop {\rm Sup}_{t\in I}}\ h(t)^{-1} \left ( \parallel \nabla B'(t)\parallel_2 \ \vee \ \parallel \partial_t B'(t)\parallel_2\right )
\eeq 

\noi where $J = [t, \infty ) \cap I$. We have already estimated $N'_i$,
$0 \leq i \leq 2$, in the proof of Proposition 2.2. We next estimate
$\nabla v'$, starting from the equation
\beq
\label{2.51e}
i \partial_t \nabla v' = - (1/2) \Delta \nabla v' + A \nabla v' + (\nabla A) v' + B \nabla u_a + (\nabla B) u_a - \nabla R_1 \ .
\eeq

\noi We first estimate $N'_3$. From (\ref{2.6e}) we obtain
\bea
\label{2.52e}
&&\parallel \nabla v'(t) \parallel_2^2\   \leq  \ \parallel \ \parallel\nabla v' \parallel_2 \Big ( \parallel \nabla A_a\parallel_{\infty}\ \parallel v'\parallel_2 \ + \ \parallel B \parallel_4 \ \parallel \nabla u_a\parallel_4\ +   \nn \\
&&+\ \parallel\nabla B\parallel_2 \ \parallel u_a \parallel_{\infty}\ + \  \parallel\nabla R_1  \parallel_2 \Big ) \ + \  \parallel \nabla v' \parallel_4 \  \parallel \nabla B \parallel_2 \ \parallel v' \parallel_4 \ ; L^1(J) \parallel \nn \\
&&\leq C \left ( N'_3 \left ( a\ N'_0 + c N_2 + c N_5\ t^{-1/2} + r_1\right ) + N'_4\ N_5\ N'_1 \ t^{-1/8}\ \overline{h}(t) \right ) h(t)^2
\eea
\noi by Lemma 2.2, so that
\beq
\label{2.53e}
N'_3 \leq C_3 \left ( a\ N'_0 + c \ N_2 + c\ N_5 \ T^{-1/2} + r_1 + \left ( N'_1\ N'_4 \ N_5 \ T^{-1/8} \overline{h}(T)\right )^{1/2}\right ) \ .
\eeq 

\noi We next estimate $N'_4$. From Lemma 2.1, we obtain

$$\parallel \nabla v';L^{8/3}(J,L^4) \parallel\   \leq  \ C \Big ( \parallel A_a \nabla v' + (\nabla A_a) v' + B \nabla u_a + (\nabla B)u_a - \nabla R_1\ ;L^1(J, L^2) \parallel$$
$$+ \ \parallel B\nabla v' + (\nabla B)v'; L^{8/5} (J, L^{4/3} )\parallel\Big )$$
$$\leq\ C \Big (  \parallel \ \parallel A_a \parallel_{\infty}\ \parallel \nabla v'\parallel_2 \ + \ \parallel \nabla A_a \parallel_{\infty} \ \parallel v' \parallel_2\ + \ \parallel B \parallel_4\ \parallel\nabla u_a\parallel_4$$
$$+ \ \parallel\nabla B\parallel_2 \ \parallel u_a \parallel_{\infty} \ +\  \parallel\nabla R_1  \parallel_2\ ; L^1(J)\parallel$$
\beq
\label{2.54e}
+\ \parallel \ \parallel B \parallel_4\ \parallel \nabla v' \parallel_2\ + \ \parallel \nabla B \parallel_2 \ \parallel v' \parallel_4\ ; L^{8/5}(J) \parallel \Big )
\eeq

\noi so that by Lemma 2.2
\bea
\label{2.55e}
&&N'_4 \leq C_4 \Big  ( a\left (  N'_3 + N'_0 \right )  + c \left ( N_2 + N_5\ T^{-1/2}\right ) + r_1 + N_2\ N'_3 \ \overline{h}(T)\nn \\
&&+ N_5\ N'_1 \ T^{-1/8} \overline{h}(T)\Big ) \ .
\eea

\noi We finally estimate $N'_5$. From (\ref{2.14e}) we obtain
\bea
\label{2.56e}
&&\parallel \nabla B(t') \parallel_2\ \vee \ \parallel \partial_t B(t') \parallel_2\ \leq \ \parallel |v|^2 + 2 \ {\rm Re} \ \overline{u}_a\ v \ + R_2; L^1(J,L^2) \parallel\nn \\
&&\leq \ \parallel \ \parallel  v\parallel _4^2\ + \ 2 \parallel  u_a \parallel_{\infty} \ \parallel  v \parallel _2 \ + \ \parallel  R_2 \parallel _2\ ;L^1(J) \parallel \nn \\
&&\leq C_5 \left ( (cN_0 + r_2)t^{-1/2} + N_1^2\ t^{1/4} \ h(t) \right ) h(t)
\eea

\noi by Lemma 2.2 so that
\beq
\label{2.57e}
N'_5 \leq C_5 \left ( (c\ N_0  + r_2)  T^{-1/2} + N^2_1\ T^{-1/8} \overline{h}(T)\right ) \ .
\eeq

We next take the limit $t_0 \to \infty$ in the same way as in Proposition 2.2. From now on we take $I = [T, \infty )$. \par

From the previous estimates it follows that the map $\phi$ defined in
Proposition 2.2, when restricted to $X_1(I)$ and more precisely to the
subset ${\cal R}_1$ of $X_1(I)$ defined by (\ref{2.19e})-(\ref{2.21e})
and (\ref{2.45e})-(\ref{2.47e}) satisfies the estimates (\ref{2.53e})
(\ref{2.55e}) (\ref{2.57e}) in addition to the previous estimates
(\ref{2.26e}) (\ref{2.29e}) (\ref{2.31e}) for $(v', B') = \phi ((v,
B))$ with $(v, B) \in {\cal R}_1$. We next show that for a suitable
choice of the $N_i$, $0 \leq i \leq 5$, and for $T$ sufficiently large,
$\phi$ maps ${\cal R}_1$ into itself. We have already chosen $N_0$,
$N_1$, $N_2$. We now have to choose $N_3$, $N_4$, $N_5$ so as to ensure
that the RHS of (\ref{2.53e}) (\ref{2.55e}) (\ref{2.57e}) do not exceed
$N_3$, $N_4$ and $N_5$ respectively. We choose 
\beq
\label{2.58e}
\left \{ \begin{array}{l} N_3 = C_3 \left ( a\ N_0 + c \ N_2 + r_1 + 1 \right )\\
N_4 = C_4 \left ( a \left ( N_3 + N_0 \right ) + c\ N_2 + r_1 + 1 \right ) \\
N_5 = C_5 \end{array} \right . 
\eeq

\noi and we take $T$ sufficiently large so that the terms not
considered in the RHS of (\ref{2.53e}) (\ref{2.55e}) (\ref{2.57e}) do
not exceed 1. Since ${\cal R}_1$ is closed in the topology of $X(I)$,
$\phi$ has a fixed point in ${\cal R}_1$ by the
contraction argument of Proposition 2.2. Finally the estimate (2.44e) follows from (\ref{2.56e}).\par \nobreak \hfill $\sq$ \par

We finally turn to the construction of solutions of the system
(\ref{1.3e}) at the level of regularity of $H^2$ for $v$. The result
can be stated as follows.\\

\noi {\bf Proposition 2.4.} {\it Let $h$ be defined as previously with
$\lambda = 3/8$ and let $X_2(\cdot )$ be defined by (\ref{1.16e}). Let
$u_a$, $A_a$, $R_1$, $R_2$ satisfy the conditions
(\ref{2.15e})-(\ref{2.18e}) and in addition 
\beq
\label{2.59e}
\parallel u_a(t) \parallel_{\infty} \ \leq  c\ t^{-3/2} \qquad , \quad \parallel \partial_t u_a (t) \parallel_4\ \leq c\ t^{-3/4} \ ,
\eeq
\beq
\label{2.60e}
\parallel \partial_t A_a \parallel_{\infty} \ \leq a \ t^{-1} \ ,
\eeq
\beq
\label{2.61e}
\parallel \partial_t R_1;L^1([t, \infty ), L^2 ) \parallel\ \leq r_1\ h(t) \ ,
\eeq
$$\parallel R_2;L^1 ([t, \infty ), L^2) \parallel\ \leq r_2\ t^{-1/2}\ h(t) \eqno(2.44)\equiv(2.62)$$

\noi for some constants $c_4$, $c$, $a$, $r_1$, $r_2$ with $c_4$
sufficiently small and for all $t \geq 1$. Then there exists $T$, $1
\leq T < \infty$ and there exists a unique solution $(v, B)$ of the
system (\ref{1.3e}) in the space $X_2([T, \infty ))$. Furthermore $B$ satisfies the estimate (2.44e). If in addition
$R_1$ satisfies the estimate
$$\parallel R_1;L^{8/3} ([t, \infty ), L^4) \parallel\ \leq r_1\ h(t) \eqno(2.63)$$

\noi for all $t \geq 1$, then $\Delta v \in L^{8/3}([T, \infty ),L^4)$ and $v$ satisfies the estimate 
$$\parallel \Delta v; L^{8/3}([t, \infty ),L^4) \parallel\ \leq\ C\ h(t)
\eqno(2.64)$$

\noi for some constant $C$ and for all $t \geq T$.}\\

\noi {\bf Proof.} The proof is essentially the same as that of Proposition
2.3 with $\nabla v$ and $\nabla v'$ replaced everywhere by $\partial_t
v$ and $\partial_t v'$ and with additional estimates of $\Delta v'$. The
existence of $v'_{t_0}$ with the required properties follows from
Proposition 2.1, part (3) and the subset ${\cal R}_2$ of $X_2(I)$
invariant under $\phi$ is now defined by the conditions
(\ref{2.19e})-(\ref{2.21e}) and in addition
$$\parallel \partial_t v(t) \parallel_2 \ \leq N_3\ h(t)
\eqno(2.65)$$
$$\parallel \partial_t v ; L^{8/3}([t, \infty ), L^4)\parallel \ \leq \  N_4\ h(t) 
\eqno(2.66)$$
$$\parallel \nabla B(t)\parallel_2 \ \vee \ \parallel \partial_t B(t)\parallel_2\ \leq \  N_5\ h(t) \eqno(2.47)\equiv(2.67)$$
$$\parallel \Delta v(t) \parallel_2\ \leq N_6\ h(t) \eqno(2.68)$$

\noi for all $t \in I$. The estimates associated with (2.65) (2.66)
(2.67) are again (\ref{2.53e}) (\ref{2.55e}) (\ref{2.57e}) except for the fact that we have to estimate in addition
$$\parallel \partial_t v'(t_0)\parallel_2 \ = \  \parallel (B u_a - R_1)(t_0)\parallel_{2} \ .\eqno(2.69)$$

\noi For that purpose, we need pointwise estimates in time of $R_1$ and $B$. From (\ref{2.61e}) it follows that $R_1 \in {\cal C}([1
, \infty ),L^2)$ and that 
$$\parallel R_1(t) \parallel_2 \ \leq \ \parallel \partial_t R_1;L^1([t, \infty ),L^2) \parallel\ \leq r_1\ h(t)\eqno(2.70)$$

\noi for all $t\geq 1$, while by the definition of ${\cal R}_2$ and by Lemma 2.2 
$$\parallel B(t) \parallel_3^3 \ \leq \ 3 \parallel \ \parallel B
\parallel_4^2 \ \parallel \partial_t B \parallel_2\ ; L^1(J) \parallel\
\leq C\ N_2^2 \ N_5\ t^{1/2}\ h(t)^3 \eqno(2.71)$$

\noi since $\parallel B(t)\parallel_3 \ \to 0$ as $t \to \infty$, which
can be proved by using a finite time version of (2.71) together  with
the fact that $\parallel B(t)\parallel_6 \ \to 0$ as $t \to \infty$ by
(\ref{2.47e}). Therefore 
$$\parallel B(t) \parallel_3\ \leq \widetilde{N}_2\ t^{1/6}\ h(t) \eqno(2.72)$$

\noi for all $t \in I$, with $\widetilde{N}_2 = C N_2^{2/3} N_5^{1/3}$. \par

We then estimate
$$\parallel Bu_a \parallel_2\ \leq \ \parallel B\parallel_3\ \parallel u_a \parallel_6\ \leq c\ \widetilde{N}_2 \ t^{-5/6}\ h(t) \ .\eqno(2.73)$$

\noi From (2.69) (2.70) (2.73) and the preceding remarks, it follows that $N'_3$ now defined by 
$$N'_3 = \ \mathrel{\mathop {\rm Sup}_{t \in I }}\ h(t)^{-1} \parallel \partial_t v'(t) \parallel_2$$

\noi satisfies an estimate obtained from (\ref{2.53e}) by adding an extra term $c\widetilde{N}_2T^{-5/6}$. \par

We have in addition to estimate $\Delta v'$. From (\ref{1.5e}) we obtain
$$\parallel \Delta v'\parallel_r\   \leq  \ 2 \left ( \parallel
\partial_t v'\parallel_r \ + \  \parallel A_a \parallel_{\infty} \
\parallel v'\parallel_{r}\ + \ \parallel Bv'\parallel_r \ + \ \parallel
B u_a\parallel_r \ + \ \parallel R_1\parallel_r\right ) \ .\eqno(2.74)$$

For $r = 2$, we estimate
\begin{eqnarray*}
\parallel Bv' \parallel_2 &\leq&  C\parallel B\parallel_3\ \parallel v' \parallel_2^{1/2}\ \parallel \Delta v'\parallel_2^{1/2}\nn \\
&\leq & (1/4) \parallel \Delta v'\parallel_2 \ + C\ \widetilde{N}_2^2\ N'_0\ t^{1/3}\ h(t)^3
\end{eqnarray*}

\noi so that from (2.74) with $r = 2$ and from (2.70) (2.73)
$$\parallel \Delta v' \parallel_2\ \leq 4 \left ( N'_3 + a\ t^{-1}\ N'_0 + r_1 + c \ \widetilde{N}_2\ t^{-5/6} + C \ \widetilde{N}_2^2\ N'_0\ t^{1/3} \ h(t)^2 \right ) h(t)\eqno(2.75)$$

\noi which suffices for the needs of the proof. \par

If $R_1$ satisfies (2.63), we can in addition derive (2.64) for $v'$. We estimate
$$\begin{array}{ll} \parallel Bu_a;L^{8/3}(J,L^4) \parallel\ &\leq \ c \parallel B;L^4(J,L^4)\parallel \ \parallel t^{-3/2} \parallel_8 \\ \\ &\leq c\ N_2\ t^{-11/8}\ h(t)\end{array} \eqno(2.76)$$

$$\parallel Bv';L^{8/3}(J,L^4) \parallel\ \leq \ \parallel B;L^4(J,L^4)\parallel \ \parallel v';L^8(J,L^{\infty}) \parallel$$
$$\leq \ C \parallel B;L^{4}(J,L^4) \parallel\ \parallel v';L^8(J,L^2)\parallel^{1/24} \ \parallel v';L^{8/3}(J,L^4) \parallel^{1/3} \ \parallel \Delta v';L^{\infty}(J,L^2) \parallel^{5/8}\ . \eqno(2.77)$$

\noi From (2.74) with $r=4$ and from (2.76) (2.77), we obtain
$$\parallel \Delta v';L^{8/3}(J,L^4) \parallel \ \leq 2 \Big ( N'_4 + a\ t^{-1}\ N'_1 + r_1 + c\ N_2 \ t^{-11/8}$$
$$+ C\ N_2\ N'^{1/24}_0\ N'^{1/3}_1\ N'_6\ h(t)^2 \Big ) h(t) \eqno(2.78)$$

\noi where in the same way as before
$$N'_6 = \ \mathrel{\mathop {\rm Sup}_{t \in I }}\ h(t)^{-1} \parallel \Delta v'(t) \parallel_2\ .$$

\noi and $N'_6$ is estimated by (2.75). \par \nobreak \hfill $\sq$ \par

\noi {\bf Remark 2.2.} There is some flexibility in the choice of the
function spaces used here. For instance we have included the Strichartz
norms $L^q(L^r)$ in the restricted range $2 \leq r \leq 4$. One could
equally well use the full range $2 \leq r\leq 6$. Conversely one could
omit the $L^2$ norm of $\partial_t B$ in the $H^1$ theory of
Proposition 2.3 and/or the $L^2$ norm of $\nabla B$ in the $H^2$ theory
of Proposition 2.4.

\mysection{Remainder estimates and completion of the proof}
\hspace*{\parindent} In this section we first prove that the choice of the asymptotic
functions $(u_a, A_a)$ made in the introduction satisfies the
assumptions of Propositions 2.2-2.4 under suitable assumptions on the
asymptotic state $(u_+, A_+, \dot{A}_+)$. We then combine those results
with those of Section 2 to complete the proof of Proposition 1.1. \par

We first supplement the definition of $(u_a, A_a)$ with some additional
properties of a general character. In addition to the representation
(\ref{1.11e}) (\ref{1.12e}) for $A_1$, we need a representation of
$\partial_t A_1$. From (\ref{1.10e}) it follows that 
\beq
\label{3.1e}
\partial_t A_1(t) = \int_t^{\infty} dt' \cos \left ( \omega (t-t')\right ) |u_a(t')|^2
\eeq

\noi so that upon substitution of (\ref{1.6e}) we obtain
\beq
\label{3.2e}
\partial_t A_1(t) = t^{-2} \ D_0(t) \widetilde{\widetilde{A}}_1
\eeq

\noi where
\beq
\label{3.3e}
\widetilde{\widetilde{A}}_1 = \int_1^{\infty} d\nu \ \nu^{-3} \cos (\omega (\nu - 1))D_0(\nu ) |w_+|^2 \ .
\eeq

\noi On the other hand from (\ref{1.11e})
\beq
\label{3.4e}
\nabla A_1(t) = t^{-2} \ D_0 (t) \nabla \widetilde{A}_1 \ .
\eeq

\noi We shall also need the commutation relations
\beq
\label{3.5e}
\nabla MD = MD \left ( ix + t^{-1} \nabla \right )
\eeq
\beq
\label{3.6e}
i \partial_t MD = MD \left ( (1/2) x^2 + i \partial_t - it^{-1} (x\cdot \nabla + 3/2)\right )
\eeq
\beq
\label{3.7e}
\left ( i \partial_t + (1/2) \Delta \right ) MD = MD \left ( i \partial_t + (2t^2)^{-1}\Delta \right ) \ .
\eeq

\noi From (\ref{3.5e}) (\ref{3.6e}) and (\ref{1.13e}), it follows that
\beq
\label{3.8e}
\nabla u_a = MD \exp (- i \varphi ) \left ( ixw_+ + t^{-1} \nabla w_+ - i t^{-1} \ell n \ t \ (\nabla \widetilde{A}_1) w_+\right )
\eeq
\bea
\label{3.9e}
&&\partial_t u_a = MD \exp (- i \varphi ) \Big ( (1/2) x^2 w_+ - i t^{-1} (x \cdot \nabla + 3/2) w_+\nn \\
&&+ t^{-1} \ \widetilde{A}_1 \ w_+ - t^{-1} \ell n \ t \ (x\cdot \nabla \widetilde{A}_1)w_+ \Big ) \ .
\eea
\noi We now consider the remainder $R_1$ defined by (\ref{1.4e}). (We recall that the choice (\ref{1.9e}) (\ref{1.10e}) yields $R_2 = 0)$. From (\ref{3.7e}) (\ref{1.13e}) we obtain
\beq
\label{3.10e}
R_1 = MD \exp (- i \varphi ) \widetilde{R}_1 - A_0 \ u_a
\eeq

\noi where
\beq
\label{3.11e}
\widetilde{R}_1 = (2t^2)^{-1} \left ( \Delta w_+ - i \ \ell n \ t \left ( 2(\nabla \widetilde{A}_1)\cdot \nabla w_+ + (\Delta \widetilde{A}_1) w_+\right ) - (\ell n\ t)^2 \ |\nabla \widetilde{A}_1|^2 w_+ \right ) \ . 
\eeq

\noi From (\ref{3.5e}) (\ref{3.6e}) (\ref{3.10e}) and (\ref{1.13e}) we obtain
\beq
\label{3.12e}
\nabla R_1 = MD \exp (-i \varphi ) \left ( ix + t^{-1} \nabla - i t^{-1} \ell n\ t \ \nabla \widetilde{A}_1\right ) \widetilde{R}_1 - \nabla (A_0 u_a)
\eeq
\bea
\label{3.13e}
&&i \partial_t R_1 = MD \exp (-i \varphi ) \Big ( (1/2) x^2 + i \partial_t - i t^{-1} (x \cdot \nabla + 3/2) \nn \\
&&+ \ t^{-1} \widetilde{A}_1 - t^{-1} \ell n\ t  \left ( x\cdot \nabla \widetilde{A}_1\right ) \Big ) \widetilde{R}_1 - i \partial_t (A_0u_a)\ .
\eea

\noi Finally, since $\widetilde{A}_1$ is independent of $t$, $\partial_t \widetilde{R}_1$ takes the explicit form 
\bea
\label{3.14e}
&&\partial_t \widetilde{R}_1 = t^{-3} \Big ( - \Delta w_+ + i (\ell n\ t - 1/2) \left ( 2 (\nabla \widetilde{A}_1)\cdot \nabla w_+ + (\Delta \widetilde{A}_1)w_+\right )\nn \\
&& + \ell n \ t (\ell n\ t -1) |\nabla \widetilde{A}_1|^2 w_+ \Big ) \ .
\eea

In order to ensure the assumptions of Propositions 2.2-2.4 on $u_a$, $A_a$, $R_1$, we shall use a number of general norm estimates. We first consider $u_a$. From (\ref{1.6e}) (\ref{3.8e}) (\ref{3.9e}) we obtain 
\beq
\label{3.15e}
\parallel u_a \parallel_r \ \leq t^{-\delta (r)} \parallel w_+\parallel_r\ , 
\eeq
\beq
\label{3.16e}
\parallel \nabla u_a \parallel_r \ \leq t^{-\delta (r)} \left ( \parallel x w_+\parallel_r\ + \ t^{-1} \parallel \nabla w_+ \parallel_r\ + \ t^{-1}\ell n\ t \parallel \nabla \widetilde{A}_1 \parallel_{\infty}\ \parallel w_+ \parallel_r\right ) \ ,
\eeq
\bea
\label{3.17e}
&&\parallel \partial_t u_a \parallel_r \ \leq t^{-\delta (r)} \Big ( (1/2) \parallel x^2 w_+\parallel_r\ + \ t^{-1} \parallel x\cdot \nabla w_+ \parallel_r\nn \\
&&+ \ t^{-1}\left ( \parallel\widetilde{A}_1 \parallel_{\infty}\ + (3/2)\right ) \parallel w_+ \parallel_r \ + t^{-1}\ell n \ t \parallel \nabla \widetilde{A}_1 \parallel_{\infty}\ \parallel xw_+ \parallel_r\Big )\Big )
\eea

\noi where $\delta (r) = 3/2 - 3/r$. \par

We next turn to $A_1$. From (\ref{1.11e}) (\ref{3.4e}) (\ref{3.2e}) we obtain
\beq
\label{3.18e}
\parallel A_1 \parallel_{\infty} \ = \ t^{-1} \parallel \widetilde{A}_1\parallel_{\infty} \ ,
\eeq
\beq
\label{3.19e}
\parallel \nabla A_1 \parallel_{\infty} \ = \ t^{-2} \parallel \nabla \widetilde{A}_1\parallel_{\infty} \ ,
\eeq
\beq
\label{3.20e}
\parallel \partial_t A_1 \parallel_{\infty} \ = \ t^{-2} \parallel \widetilde{\widetilde{A}}_1\parallel_{\infty} \ .
\eeq

\noi The $L^{\infty}$ estimates of $\widetilde{A}_1$ and $\widetilde{\widetilde{A}}_1$ will be obtained through Sobolev inequalities from the $L^2$ estimates
\bea
\label{3.21e}
\parallel \omega^{m+1} \widetilde{A}_1\parallel_2\ \vee \ \parallel\omega^m \widetilde{\widetilde{A}}_1\parallel_2 &\leq & \int_1^{\infty} d\nu \ \nu^{-3/2-m} \parallel \omega^m |w_+|^2\parallel_2\nn \\
&=& (m+1/2)^{-1} \parallel \omega^m |w_+|^2\parallel_2
\eea

\noi which follow readily from (\ref{1.12e}) (\ref{3.3e}). Finally we shall estimate $R_1$ and its derivatives as follows~:
\beq
\label{3.22e}
\parallel R_1 \parallel _2\ \leq \ \parallel \widetilde{R}_1\parallel_2 \ + \ t^{-3/2} \parallel A_0 \parallel _2 \ \parallel  w_+ \parallel _{\infty}
\eeq

\noi where we have used (\ref{3.15e}) and where
\bea
\label{3.23e}
&&\parallel \widetilde{R}_1\parallel_2\ \leq (2t^2)^{-1} \Big ( \parallel \Delta w_+ \parallel_{2} \ + \ \ell n\ t \ \Big ( 2\parallel \nabla \widetilde{A}_1 \parallel_{6}\ \parallel \nabla w_+ \parallel_{3}\nn \\
&&+\ \parallel \Delta \widetilde{A}_1 \parallel_{2}\ \parallel w_+ \parallel_{\infty}\Big ) + (\ell n\ t )^2 \parallel \nabla \widetilde{A}_1 \parallel_{6}^2 \ \parallel \nabla w_+ \parallel_{6} \Big )\ ,
\eea
\bea
\label{3.24e}
&&\parallel \nabla R_1 \parallel_{2}\ \leq \ \parallel x \widetilde{R}_1 \parallel_{2}\ + \ t^{-1} \parallel \nabla \widetilde{R}_1 \parallel_{2}\ + \ t^{-1} \ell n\ t \parallel \nabla \widetilde{A}_1 \parallel_{\infty}\ \parallel \widetilde{R}_1 \parallel_{2}\nn \\
&&+ \ t^{-3/2} \parallel \nabla A_0 \parallel_{2}\ \parallel w_+ \parallel_{\infty}\ + \ \parallel A_0\parallel_{2}\ \parallel \nabla u_a \parallel_{\infty}\ ,
\eea
\bea
\label{3.25e}
&&\parallel \partial_t R_1 \parallel_{2}\ \leq (1/2) \parallel x^2 \widetilde{R}_1 \parallel_{2}\ + \ \parallel \partial_t \widetilde{R}_1 \parallel_{2}\ + t^{-1} \parallel x\cdot \nabla \widetilde{R}_1 \parallel_{2}\nn \\
&&+ \ t^{-1}\left ( \parallel \widetilde{A}_1 \parallel_{\infty}\ + 3/2 \right ) \parallel \widetilde{R}_1 \parallel_{2}\ + \ t^{-1}\ell n \ t \parallel \nabla \widetilde{A}_1 \parallel_{\infty}\ \parallel x \widetilde{R}_1 \parallel_{2}\nn \\
&&+ \ t^{-3/2}\parallel \partial_t A_0 \parallel_{2}\ \parallel w_+ \parallel_{\infty}\ + \ \parallel A_0 \parallel_{2}\ \parallel \partial_t u_a \parallel_{\infty}\ ,
\eea

\noi where $\nabla u_a$ and $\partial_t u_a$ are estimated in $L^{\infty}$ by (\ref{3.16e}) (\ref{3.17e}) with $r=\infty$. \par

In order to estimate $A_0$, we need some general estimates of solutions of the free wave equation.\\

\noi {\bf Lemma 3.1.} {\it Let $A_0$ be defined by (\ref{1.9e}). Let $k \geq 0$ be an integer. Let $A_+$ and $\dot{A}_+$ satisfy the conditions
\beq
\label{3.26e}
A_+, \omega^{-1}\dot{A}_+ \in H^k \qquad , \quad \nabla^2 A_+, \nabla \dot{A}_+ \in W_1^k \ .
\eeq
\noi Then $A_0$ satisfies estimates
\beq
\label{3.27e}
\left \{ \begin{array}{l}  \parallel A_0(t);W_r^k \parallel \ \leq a_0\ t^{-1+2/r} \ ,\\ \\ \parallel \partial_t A_0(t);W_r^{k-1}\parallel\ \leq a_0\ t^{-1+2/r} \quad \hbox{\it for $k \geq 1$}\ .\end{array} \right .
\eeq

\noi for $2 \leq r \leq \infty$ and for all $t\in {I\hskip-1truemm R}$, where $a_0$
depends on $A_+$, $\dot{A}_+$ through the norms associated with
(\ref{3.26e}).}\\

A proof can be found in \cite{17r}.\par

We are now in a position to derive the estimates required in
Propositions 2.2-2.4. We recall that $w_+ = Fu_+$ and that $\delta (r)
= 3/2-3/r$.\\

\noi {\bf Proposition 3.1.} \par
{\it (1) Let $u_+ \in H^{0,2}$ and let $(A_+, \dot{A}_+)$ satisfy (\ref{3.26e}) with $k = 0$. Then the following estimates hold
\beq
\label{3.28e}
\parallel u_a\parallel_r\ \leq t^{-\delta (r)} \parallel Fu_+ \parallel_r \qquad \hbox{\it for $2 \leq r \leq \infty$} \ ,
\eeq
\beq
\label{3.29e}
\parallel A_a(t)\parallel_{\infty} \ \leq a\ t^{-1} \ ,
\eeq
\beq
\label{3.30e}
\parallel R_1(t)\parallel_2 \ \leq r_1\ t^{-3/2}\ .
\eeq

(2) Let $u_+ \in H^{0,3} \cap H^{1,2}$ and let $(A_+, \dot{A}_+)$ satisfy
(\ref{3.26e}) with $k = 1$. Then the estimates
(\ref{3.28e})-(\ref{3.30e}) hold and in addition
\beq
\label{3.31e}
\parallel \nabla u_a\parallel_r\ \leq c\ t^{-\delta (r)} \qquad \hbox{\it for $2 \leq r \leq \infty$}\ ,
\eeq
\beq
\label{3.32e}
\parallel \nabla A_a(t)\parallel_{\infty} \ \leq a\ t^{-1} \ ,
\eeq
\beq
\label{3.33e}
\parallel \nabla R_1(t)\parallel_2 \ \leq r_1\ t^{-3/2}\ .
\eeq

(3) Let $u_+ \in H^{1,3}\cap H^{2,2}$ and let $(A_+, \dot{A}_+)$
satisfy (\ref{3.26e}) with $k = 1$. Then the estimates
(\ref{3.28e})-(\ref{3.33e}) hold and in addition}
\beq
\label{3.34e}
\parallel \partial_t u_a\parallel_r\ \leq c\ t^{-\delta (r)} \qquad \hbox{\it for $2 \leq r \leq \infty$} \ ,
\eeq
\beq
\label{3.35e}
\parallel \partial_t A_a(t)\parallel_{\infty} \ \leq a\ t^{-1} \ ,
\eeq
\beq
\label{3.36e}
\parallel \partial_t R_1(t)\parallel_2 \ \leq r_1\ t^{-3/2}\ .
\eeq

\noi {\bf Proof.} \underbar{Part (1)}. The assumption on $u_+$ is
equivalent to $w_+ \in H^2$, which by (\ref{3.21e}) implies that
$\nabla \widetilde{A}_1 \in H^2$. The estimate (\ref{3.28e}) is a
rewriting of (\ref{3.15e}) and is ensured by the fact that $H^2 \subset
L^{\infty}$. The estimate (\ref{3.29e}) follows from (\ref{3.27e}) as
regards the $A_0$ part and from (\ref{3.18e}) and the previous remarks
as regards the $A_1$ part. Finally (\ref{3.30e}) follows from
(\ref{3.22e}) (\ref{3.23e}) (\ref{3.27e}) and Sobolev inequalities. \\

\underbar{Part (2)}. The assumption on $u_+$ is equivalent to $w_+ \in
H^3 \cap H^{2,1}$ which by (\ref{3.21e}) implies that $\nabla
\widetilde{A}_1 \in H^3$. Then (\ref{3.31e}) follows from (\ref{3.16e})
while (\ref{3.32e}) follows from (\ref{3.27e}) as regards the $A_0$
part and from (\ref{3.19e}) and the previous remarks as regards the
$A_1$ part. Finally (\ref{3.33e}) follows from (\ref{3.24e}). In fact,
in addition to terms previously estimated, we have to estimate
$\parallel x \widetilde{R}_1\parallel_2$ and $\parallel \nabla \widetilde{R}_1\parallel_2$. The estimate of
$x\widetilde{R}_1$ is obtained from (\ref{3.23e}) by replacing $\Delta
w_+$, $\nabla w_+$ and $w_+$ by $x \Delta w_+$, $x\nabla w_+$ and
$xw_+$ respectively. The estimate of $\nabla \widetilde{R}_1$ is
obtained from (\ref{3.23e}) by distributing $\nabla$ among
$\widetilde{A}_1$ and $w_+$, thereby generating norms which are
controlled by the assumption $w_+ \in H^3$.\\

\underbar{Part (3)}. The assumption on $u_+$ is equivalent to $w_+ \in
H^{3,1} \cap H^{2,2}$ which by (\ref{3.21e}) implies that $\nabla
\widetilde{A}_1 \in H^3$ and $\widetilde{\widetilde{A}}_1 \in H^2$.
Then (\ref{3.34e}) follows from (\ref{3.17e}) while (\ref{3.35e})
follows from (\ref{3.27e}) as regards the $A_0$ part and from
(\ref{3.20e}) and the previous remarks as regards the $A_1$ part.
Finally (\ref{3.36e}) follows from (\ref{3.25e}). In fact, in addition
to terms previously estimated, we have to estimate $\parallel
x^2\widetilde{R}_1\parallel_2$, $\parallel x \cdot \nabla
\widetilde{R}_1\parallel_2$ and $\parallel \partial_t
\widetilde{R}_1\parallel_2$. In the same way as in the proof of Part
(2), the first two estimates are obtained from that of $\parallel
\widetilde{R}_1\parallel_2$ by absorbing $x^2$ or $x$ by $w$ and
distributing the gradient among $w_+$ and $\widetilde{A}_1$, while $\parallel \partial_t
\widetilde{R}_1\parallel_2$ is estimated in the same way as $\parallel 
\widetilde{R}_1\parallel_2$ from the explicit expression (\ref{3.14e}).\par \nobreak \hfill $\sq$ \par

We can now complete the proof of Proposition 1.1. \\

\noi {\bf Proof of Proposition 1.1}. From Parts (1), (2) and (3) of
Proposition 3.1, together with the fact that $R_2 = 0$, it follows that
the assumptions of Propositions 2.2, 2.3 and 2.4 respectively are
satisfied with $h(t) = t^{-1/2}$ and $c_4 = \parallel w_+\parallel_4$.
In particular (\ref{3.30e}) (\ref{3.33e}) imply (2.63) since
\beq
\label{3.37e}
\parallel R_1(t)\parallel_4\ \leq \ C\parallel R_1(t);H^1\parallel\ \leq C\ t^{-3/2}
\eeq

\noi so that
\beq
\label{3.38e}
\parallel R_1;L^{8/3}([t, \infty ),L^4\parallel\ \leq \ C\ t^{-9/8} \ .
\eeq

\noi The estimate (\ref{1.17e}) follows from (2.44e) with $h(t) = t^{-1/2}$.
\par\nobreak \hfill $\sq$ \par

\noi {\bf Remark 3.1.} The $t^{-3/2}$ decay of $R_1$ comes from the
free wave term $A_0u_a$. That term could be partly cancelled by the
correcting term used in \cite{15r}, thereby producing a $t^{-2}(\ell n\
t)^2$ decay of $R_1$ allowing for $h(t) = t^{-1}(\ell n\ t)^2$. \\

\noi {\bf Remark 3.2.} The regularity assumptions on $u_+$ or $w_+$ are
dictated by the term $\Delta w_+$ in $\widetilde{R}_1$. They could be
somewhat weakened by eliminating that term through the choice $$w(t) =
U(1/t)^* w_+$$

\noi but that choice would either generate a more complicated and less
explicit $\varphi$ or produce a non vanishing $R_2$.

\end{document}